\documentstyle[12pt]{article}
\input amssym.def
\input amssym
\def\qed{\nopagebreak\par\noindent\nopagebreak$\blacksquare$\par}
\def\reals{{\Bbb R}}
\def\naturals{{\Bbb N}}
\def\proof{{\par\noindent{\bf Proof: }}}
\def\base{{\cal B}}
\newtheorem{theorem}{Theorem}
\newtheorem{lemma}[theorem]{Lemma}

\begin{document}

\bigskip
\begin{center}
{\Large  Souslin's Hypothesis and Convergence in Category}
\end{center}
\bigskip

\begin{center}
 by Arnold W. Miller\footnote{
I want to thank Krzysztof Ciesielski for many
helpful conversations.

The results presented in this paper were obtained
during the Joint US--Polish Workshop in Real Analysis,
{\L}{\'o}d{\'z}, Poland, July 1994.
The Workshop was partially
supported by the NSF grant INT--9401673.

  }\\
\end{center}

Suppose that $S\subseteq P(X)$ is a $\sigma$-field of
subsets of $X$ and $I\subseteq S$ is a $\sigma$-ideal.
If $I$ has
the countable chain condition (ccc),  i.e., every family of
disjoint sets in $S\setminus I$ is countable, then
$S/I$ is a complete boolean algebra.  A boolean algebra
is atomic iff there is an atom beneath every nonzero element.

A function $f:X\to \reals$ is $S$-measurable iff
$f^{-1}(U)\in S$ for every open set $U$.
A sequence of $S$-measurable functions $f_n:X\to \reals$
converges $I$-a.e. to a function $f$ iff there exists $A\in I$ such
that $f_n(x)\to f(x)$ for all $x\in (X\setminus A)$.
If $(X,S,\mu)$ is a finite measure space, then a sequence of
measurable functions $f_n:X\to \reals$ converges in measure
to a function $f$ iff for any $\epsilon>0$ there exists $N$ such
that for any $n>N$:
$$\mu(\{x\in X: |f_n(x)-f(x)|>\epsilon\})<\epsilon.$$
In this case if $I$ is the ideal of measure zero sets, then
$f_n$ converges to $f$ in measure iff
every subsequence $\{f_n:n \in X\}$ (where $X\subseteq \naturals$)
has a subsequence $Y\subseteq X$ such that $\{f_n: n\in Y\}$
converges $I$-a.e.   This allows us to define convergence in
measure without mentioning the measure,  only  the
ideal $I$.  So in the abstract setting define $f_n$ converges
to $f$ with respect to $I$ iff every subsequence $\{f_n:n \in X\}$
has a subsequence $Y\subseteq X$ such that $\{f_n: n\in Y\}$
converges $I$-a.e. (where $X$ and $Y$ range over infinite sets of
natural numbers.)  For more background on this subject in
case $I$ is the ideal of meager sets, see
Poreda, Wagner-Bojakoska, and Wilczy\'{n}ski [PWW] and
Ciesielski, Larson, and Ostaszewski [CLO].

Marczewski [M] showed that
if $(X,S,\mu)$ is an atomic measure and $I$ the
$\mu$-null sets, then
`$I$-a.e. convergence' is the same as `convergence with respect to I'.

Gribanov [G] proved the converse,
if $(X,S,\mu)$ is a finite measure space and $I$ the
$\mu$-null sets, then
if `$I$-a.e. convergence' is the same as `convergence with respect to I'
then  $\mu$ is an atomic measure.

Souslin's Hypothesis (SH) is the statement that there are no
Souslin lines.  It is known to be independent (see Solovay and
Tennenbaum [ST]).  It was the inspiration for Martin's Axiom.

\begin{theorem} (Wagner and Wilczy\'{n}ski [WW]) \label{WW}
Assume SH.
Then for any $\sigma$-field $S$ and ccc $\sigma$-ideal $I\subseteq S$
the following are equivalent:
\begin{itemize}
  \item `$I$-a.e. convergence' is the same as `convergence with respect to I'
  for $S$-measurable sequences sequences of real-valued functions, and
  \item the complete boolean algebra $S/I$ is atomic.
\end{itemize}
\end{theorem}

At the real analysis meeting in \L\'{o}d\'{z} Poland in July 94,
Wilczy\'{n}ski asked whether or not SH is needed for the Theorem
above.  We show here that the conclusion of Theorem \ref{WW} implies
Souslin's Hypothesis.

\begin{theorem}
 Suppose SH is false (so there exists a Souslin tree).  Then
 there exists a regular topological space $X$ such that
 \begin{enumerate}
   \item $X$ has no isolated points,
   \item $X$ is ccc (every family of disjoint open sets is countable),
   \item  every nonempty open subset of $X$ is nonmeager, and
   \item  if $I$ is
 the $\sigma$-ideal of meager subsets of $X$, then
 `$I$-a.e. convergence' is the same as `convergence with respect to I'
 for any sequence of Baire measurable real-valued functions.
 \end{enumerate}
Hence if $S$ is $\sigma$-ideal of sets with the property of Baire
and $I$ the $\sigma$-ideal of meager sets, then $S/I$ is ccc and
nonatomic, but the two types of convergence are the same.
\end{theorem}
\proof
Define $(T,<)$ is an $\omega_1$-tree iff it is a partial order and
for each $s\in T$ the set $\{t\in T: t< s\}$ is well-ordered
by $<$ with some countable order type, $\alpha<\omega_1$.
We let
$$T_\alpha=\{s\in T: \{t\in T: t<s\}
\mbox{ has order type } \alpha\}.$$
Also
$$T_{<\alpha}=\bigcup\{T_\beta:\beta<\alpha\}.$$

Define $C\subseteq T$ is a chain iff for
every $s,t\in C$ either $s\leq t$ or $t\leq s$.

Define
$A\subseteq T$ is an antichain iff for any $s,t\in A$
if $s\leq t$, then  $s=t$, i.e. distinct elements are
$\leq$-incomparable.

Define $T$ is a Souslin tree iff
$T$ is an $\omega_1$ tree in which every chain and antichain
is countable.   (Note that since $T_\alpha$ is an
antichain it must be countable.)

SH is equivalent to saying there is no Souslin
tree.  Every Souslin tree contains
a normal Souslin tree, i.e., a Souslin tree $T$ such that
for every $\alpha<\beta<\omega_1$ and $s\in T_\alpha$ there
exists a $t\in T_\beta$ with $s<t$.  (Just throw out nodes
of $T$ which do not have extensions arbitrarily high in the tree.)
For more on Souslin trees see
Todor\v{c}evi\v{c} [T].

Now we are ready to define our space $X$.  Let the elements of
$X$ be maximal chains of $T$.  For each $s\in T$ let
$$C_s=\{b\in X: s\in b\}$$
and let
$$\{C_s:s\in T\}$$
be an open basis for the topology on $X$.  Note that
$C_s\cap C_t$ is either empty or equal to either
$C_s$ or $C_t$ depending on whether $s$ and $t$ are incomparable,
or $t\leq s$ or $s\leq t$, respectively.
Each $C_s$ is clopen since its complement is the union
of $C_t$ for $t$ which are incomparable to $s$.  $X$ has no isolated
points, since given any $s\in T$ there must be incomparable extensions
of $s$ (because $T$ is normal) and therefore at
least two maximal chains containing $s$, so $C_s$ is
not a singleton. Clearly $X$ has the countable chain condition.

\begin{lemma}\label{baire}
Open subsets of $X$ are nonmeager.  In fact, the intersection of
countably many open dense sets contains an open dense set.
\end{lemma}
\proof
Suppose $(U_n:n\in\omega)$ is a sequence of an open dense subsets of $X$.
Let $A_n\subseteq T$ be
an antichain which is maximal with respect to the property
that $C_s\subseteq U_n$ for each $s\in A_n$.
Since $U_n$ is open dense in $X$, $A_n$ will be a
maximal antichain in $T$.

Let
$$V_n=\bigcup\{C_s:s\in A_n\}.$$
Then $V_n\subseteq U_n$ and $V_n$ is open dense.  (It is dense, because
given any $C_t$ there exists $s\in A_n$ and $r\in T$ with
$t\leq r$ and $s\leq r$, hence $C_r\subseteq V_n\cap C_t$.)

Choose $\alpha<\omega_1$ so that for each $n\in\omega$ the
(necessarily countable) antichain $A_n\subseteq T_{<\alpha}$.
Let $$U=\bigcup\{C_s:s\in T_\alpha\}.$$  Note that since $T$ is normal
$U$ is an open dense set.  Also
$$U\subseteq \bigcap_{n<\omega} V_n \subseteq \bigcap_{n<\omega} U_n.$$
($U\subseteq V_n$ because for any $b\in U$ if $b\in C_s$ for
some $s\in T_\alpha$ there must be $t\in A_n$ comparable to it, since
$A_n$ is a maximal antichain, and since $A_n\subseteq T_{<\alpha}$,
it must be that $t<s$ and so $b\in C_t\subseteq V_n$.
\qed

\begin{lemma} \label{func}
Suppose $f:X\to\reals$ is a
real valued Baire function.  Then there exists $\alpha<\omega_1$
such that for each $s\in T_\alpha$ the function
$f$ is constant on $C_s$.
\end{lemma}
\proof
Let $\base$ be a countable open basis for $\reals$.  For each
$B\in\base$  the set
$f^{-1}(B)$ has the property of Baire (open modulo meager).
So there exists an open $U_B$ such that
$$U_B \Delta f^{-1}(B) \mbox{ is meager.}$$
By the proof of Lemma \ref{baire} we may assume that
$$U_B=\bigcup\{C_s: s\in A_B\}$$ for some countable set $A_B\subseteq T$.
By the proof of Lemma \ref{baire} there exists an $\alpha<\omega_1$
such that
\begin{itemize}
  \item each $A_B\subseteq T_{<\alpha}$ and
  \item if $U$ is the open dense set $\bigcup\{C_s:s\in T_\alpha\}$,
  then $U$ is disjoint from $U_B \Delta f^{-1}(B)$ for each $B\in\base$.
\end{itemize}
But now, $f$ is constant on each $C_s$ for $s\in T_\alpha$.  Otherwise,
suppose that $f(b)\not=f(c)$ for some $b,c\in C_s$ for
some $s\in T_\alpha$. Then suppose that $f(b)\in B$ and $f(c)\notin B$ for
for some $B\in\base$.  Because $b\in (f^{-1}(B)\cap U)$ and
$U$ is disjoint from $U_B \Delta f^{-1}(B)$, it must be that $b\in U_B$.
Hence there exists $t\in T_{<\alpha}$ such that $C_t\subseteq U_B$ and
$b\in C_t$. Since $t<s$ it must be that $c\in C_t$ and
so $c\in f^{-1}(B)$, which contradicts $f(c)\notin B$.
\qed

Steprans [S] shows that every continuous function on
a Souslin tree takes on only countably many values.

\begin{lemma} \label{func2}
Suppose $\{f_n:X\to\reals:n\in\omega\}$ is a countable set of
real valued Baire functions.  Then there exists $\alpha<\omega_1$
such that for each $s\in T_\alpha$ and $n<\omega$ the function
$f_n$ is constant on $C_s$.
\end{lemma}
\proof
Apply Lemma \ref{func} countably many times and take the supremum
of the $\alpha_n$.
\qed

Finally, we prove the theorem.  The idea of the
proof is to use the argument of the atomic case, where the `atoms'
are supplied by Lemma \ref{func2}.
Since `$I$-a.e. convergence' always
implies `convergence with respect to $I$', it is enough to see the
converse.  So let $f_n:X\to\reals$ be Baire functions which converge
to $f:X\to\reals$ with respect to $I$, i.e. every subsequence
has a subsequence which converges on a comeager set to $f$.  By
Lemma \ref{func2} there exists $\alpha<\omega_1$
such that for each $s\in T_\alpha$ and $n<\omega$ the function
$f_n$ is constant on $C_s$.  It must be that for each fixed
$s\in T_\alpha$ the constant values of $f_n$ on $C_s$ must converge
on a comeager subset of $C_s$ to a constant value.
It follows that the sequence $f_n(x)$ converges to $f(x)$ on a comeager
subset of the dense set $\{C_s:s\in T_{\alpha}\}$.
\qed

\begin{center}
References
\end{center}

\noindent
[CLO] K.Ciesielski, L.Larson, and K.Ostaszewski,
{\bf \cal{I}-Density Continuous Functions}, Memoirs of
the American Mathematical Society, 107(1994).

\medskip\noindent
[PWW] W.Poreda, E.Wagner-Bojakoska, and W.Wilczy\'{n}ski, A category
analogue of the density topology, Fundamenta Mathematicae,
125(1985), 167-173.

\medskip\noindent
[S] J.Stepra\'{n}s, Trees and continuous mappings into the real line,
Topology and Its Applications, 12(1981), 181-185.

\medskip\noindent
[ST] R.M.Solovay, S.Tennenbaum, Iterated Cohen extensions and
Souslin's problem, Annals of Mathematics, 94(1971), 201-245.

\medskip\noindent
[T] S.Todor\v{c}evi\v{c}, Trees and linearly ordered sets, in
{\bf Handbook of Set Theoretic Topology}, North-Holland, (1984), 235-293.

\medskip\noindent
[WW] E.Wagner, W.Wilczy\'{n}ski, Convergence of sequences of
measurable functions, Acta Mathematica Academiae Scientiarum Hungaricae,
36(1980), 125-128.

\bigskip\noindent Address:
 York University,  Department of Mathematics,
 North York,  Ontario M3J 1P3, Canada.

\medskip\noindent Permanent address:
  University of Wisconsin-Madison,
    Department of Mathematics,
     Van Vleck Hall,
      480 Lincoln Drive,
       Madison, Wisconsin 53706-1388, USA).

\medskip\noindent  e-mail: miller@math.wisc.edu

\begin{center}
Nov 94
\end{center}
\end{document}